\documentclass[11pt,reqno]{amsart}
\usepackage{amstext,amsmath,amsfonts,amssymb,amsthm, latexsym, verbatim, eucal}
\newcommand{\trace}{\mathrm{tr}}
\newcommand{\ad}{\mathrm{ad}}

\newcommand{\lie}{\mathcal{L}}
\newcommand{\liet}{\lie_{\Theta}}

\newcommand{\liend}{\lie(n,d)}

\newcommand{\EEnd}{\operatorname{End}}
\newcommand{\GL}{\operatorname{GL}}

\newcommand{\Matnk}{{Mat_n(k)}}

\newcommand{\Dn}{D^{\prime}_n}

\newcommand{\Char}{\operatorname{char}}

\DeclareMathOperator{\spn}{span}
\DeclareMathOperator{\im}{im} \DeclareMathOperator{\End}{End}
\DeclareMathOperator{\Cent}{Cent} 
\newtheorem{thm}{Theorem}[section]
\newtheorem{definition}[thm]{Definition}
\newtheorem{prop}[thm]{Proposition}

\newtheorem{lem}[thm]{Lemma}
\newtheorem{rem}[thm]{Remark}

\newtheorem{ex}[thm]{Example}
\newtheorem{qu}[thm]{Question}
\begin{document}

\title{The Automorphism Group of Certain Higher Degree Forms}
\author{M. O'Ryan}
\address{Instituto de Matematica y Fisica\newline  
Universidad de Talca,\newline Casilla 721 Talca, Chile}
\email{moryan@inst-mat.utalca.cl}
\thanks{Partially supported by FONDECYT \# 1080344, \# 1040670, by Programa Reticulados y 
Ecuaciones, U. de Talca and by PBCT \#ACT 05}
\author{S. Ryom-Hansen}
\address{Instituto de Matematica y Fisica\newline  
Universidad de Talca,\newline Casilla 721 Talca, Chile}
\email{steen@inst-mat.utalca.cl}
\thanks{Partially supported by FONDECYT \# 1051024, by Programa Reticulados y Ecuaciones, U. de Talca and by PBCT \#ACT 05}

\date{}
\begin{abstract}
\noindent We consider symmetric indecomposable $d$-linear ($d>2$) spaces of dimension $n$ over an algebraically closed field $k$ of  characteristic 0, whose center (the analog of the space of symmetric matrices of a bilinear form) is cyclic, as introduced by Reichstein \cite{rich}.  The automorphism group of these spaces is determined through the action on the center and through the determination of 
the Lie algebra. Furthermore, we relate the Lie algebra to the Witt algebra.
\end{abstract}
\maketitle
\section{Introduction}
A higher degree form, that is a form of degree $ d \ge 3 $ over a field $ k $, is a map 
$\phi:V\longrightarrow k$  satisfying $\phi(\alpha v)=\alpha^d\phi(v)$  for all $\alpha\in k, \, v\in V $.
The map $\phi$ can be polarized to obtain  a symmetric $d$-linear form  $\theta:V\times\cdots\times V\longrightarrow k$. 
If $k$ has characteristic $0$ or greater than $d,$ then this construction allows to establish a correspondence between forms of degree $d$ and $d$-linear spaces $ (V, \Theta) $.
Using this correspondence  
D. Harrison  initiated an algebraic theory of higher degree forms (somehow similar to the algebraic theory of  quadratic forms). 

\medskip
Harrison, in  \cite{ha} introduced the concepts of regularity, indecomposability of a higher degree form and also defined its  center $\Cent(\Theta),$ the analog to the space of symmetric matrices with respect to a bilinear form. It turns out that the center is a commutative subalgebra of $\End(V).$  The 
automorphism group $\mathcal{O}(\Theta)$ of a regular $d$-linear spaces $(V,\Theta)$ is 
defined as the group of all linear bijections $\sigma:V\longrightarrow V$ such that $\Theta(v_1,\ldots,v_d)=\Theta(\sigma(v_1),\ldots,\sigma(v_d))$ for all $v_i\in V.$
Sometimes this group is called the {\em orthogonal group} of $(V,\Theta)$ 
in analogy with the quadratic case. The orthogonal group of a higher degree form  has been studied by several authors, for instance see \cite{cl, csw, schnei,sla,summ,suzu}.

\medskip
In this paper we shall be interested in an action $\mathcal{O}(\Theta)$ on $ \Cent(\Theta) $. In the case of forms of cyclic center, this action and its 
induced action at Lie algebra level enables us to give a description of  $\mathcal{O}(\Theta)$.

\medskip
The organization of the paper is as follows. In the next section we recall the definitions from the literature of the objects mentioned above. We also 
introduce the concept of forms with maximal center and give examples of forms satisfying this maximality property. 
We furthermore explain how an element of the center may be used to twist the Lie algebra of the form. 
In the third section  we recall the definition of cyclic forms. We show that they have maximal center and use the twisting operation to give 
a precise description of the Lie algebra. This turns out to be quotient of the positive part of the Witt algebra. 
In the last section, we combine the action of $\mathcal{O}(\Theta)$ on $ \Cent(\Theta) $ with the information on the Lie algebra gathered in the 
previous section to give a description of $\mathcal{O}(\Theta)$ in terms of a short exact sequence involving the automorphism group of 
the $ k $-algebra $ \Cent(\Theta) $. 
We  finish the paper by calculating a concrete example.

\medskip

It is a pleasure to thank E. Backelin for many useful discussions.  We would also like to thank the referee for several remarks made to the first version of this work.

\section{Forms of degree higher than 2}\label{basic}

In this section we recall the relevant definitions and basic facts from the theory of higher degree forms. 
We introduce the concept of maximal center and we explain how the center in general can 
be made into a module for the Lie algebra in a canonical way. Finally, we show how it can be used to twist 
the bracket of the Lie algebra. 
 
\medskip

Let $d$ be an integer $d\geq 2$ and $k$ a field of characteristic $0$ or greater than $d$. 

\begin{definition}\label{dlindef} 
A \emph{$d$-linear space over $k$}\/ is a pair $(V,\Theta)$
where $V$ is a finite dimensional $k$-vector space and $\Theta:
V\times\cdots\times V\longrightarrow k$ is a symmetric $d$-linear form.
\end{definition}
This means that $\Theta(v_1,\ldots,v_d)$ is $k$-linear in each of its $d$ slots
and is invariant under all permutations of those slots. 

\medskip

Under the conditions on $k$, it is known that there is a correspondence between $d$-linear spaces $(V,\Theta)$ over $k$ 
and homogeneous polynomials of degree $d$ in $n$ unknowns and coefficients in $k,$ see \cite{ha, pros}. 

\medskip

Two $d$-linear spaces $(V, \Theta)$ and $(V', \Theta')$ over $k$ are called \emph{isomorphic}, 
if there is a bijective $k$-linear map $f : V \longrightarrow V'$ such that
$\Theta'(f(v_1), \dots, f(v_d))=\Theta(v_1, \dots, v_d)$ for every $v_1, \dots,
v_d \in V$. In this case we write $(V, \Theta)\cong (V', \Theta')$.  

\medskip

If $(V,\Theta)$ is a $d$-linear space over $k$ and $K/k$ is a field extension one gets a $d$-linear space over $K$ by extension of scalars. We denote this space by  $(V_K,\Theta_K).$   

\medskip

The {\em orthogonal sum} of $d$-linear spaces 
$(V_1,\Theta_1) $ and $ (V_2,\Theta_2)$
is denoted $(V_1,\Theta_1) \perp (V_2,\Theta_2)$. It 
is the $d$-linear space on the direct sum $V_1\oplus V_2$ with
map $\Theta_1\perp\Theta_2$ given by $(\Theta_1\perp\Theta_2)(v_1+u_1,\ldots,v_d+u_d):= \Theta_1(v_1,\ldots,v_d)+\Theta_2(u_1,\ldots,u_d).$ 
A $ d$-linear space is called {\it decomposable} if it is isomorphic to the orthogonal sum of two $d$-linear spaces, otherwise 
it is called {\it indecomposable}.

\medskip

Let $ (V,\Theta) $ be a $ d$-linear space. Two subspaces $S$ and $T$ of $V$ are called \emph{orthogonal} if $\Theta(S,T,V,\dots,V) =
0$, that is $\Theta(s,t,v_3,\dots,v_d)=0$ for every $s\in S, t\in T$ and
every $v_i\in V$. For example, if $V=V_1\perp V_2$ and we consider $ V_j $ as subspaces of $ V $, 
then  $V_1$ and $V_2$ are orthogonal. Conversely, if $S, T$ are orthogonal
subspaces of $ V $ with $S+T=V$ and $S\cap T=0$, then $V\cong S\perp T$.

\medskip

\begin{definition} Let $(V,\Theta)$ be a $d$-linear space and suppose $S\subseteq V$. Then we define $S^\perp :=\{w\in V:
\Theta(w,S,V,\ldots,V)=0\}.$ We say that $(V,\Theta)$ is 1-regular (or just regular) 
if $V^\perp=0$. We say that $(V,\Theta)$ is 2-regular
if $\Theta(w,w,V,\ldots,V)=0$ implies $w=0$.

\end{definition}

\medskip
Of course $2$-regularity implies regularity and in fact 
we shall only consider spaces that are at least regular and of degree $d\geq 3$ from now on.
Harrison proved that every regular $d$-linear space $ (V, \Theta)$ can be expressed as an orthogonal sum of unique  
regular indecomposable spaces, \cite[Proposition 2.3]{ha}. Moreover, he showed that
the indecomposable components of $ (V, \Theta)$ are determined by the ``center'', which 
is defined analogously to the space of symmetric matrices for a bilinear form:

\begin{definition}
\label{centdef} Let $(V,\Theta)$ be a (regular) $d$-linear space over $k$. The \emph{center} 
is defined as 
$ \Cent_k(V,\Theta) := \{f\in \End_k(V): \Theta(f(v_1),v_2,\ldots,v_d)=
\Theta(v_1,f(v_2),\ldots,v_d) \;  for\; all\; v_i\in V\} $.
\end{definition}

Of course, by the symmetry of the form, the second slot can be replaced by the $i$'th slot for $ i \geq 2 $ in the 
above definition.  
Let us mention the following useful properties of the center of regular spaces, also established by Harrison, \cite[Section 4]{ha}. Note that $d>2$. 

\medskip

1) $\Cent_k(\Theta)$ is a commutative $k$-subalgebra of $\End_k(V)$ containing $k.$  \newline
2)  $\Cent_k((V_1,\Theta_1)\perp (V_2,\Theta_2))\cong \Cent_k(V_1,\Theta_1)\times \Cent_k(V_2,\Theta_2)$ as $k$-algebras. \newline
3)  $\Cent_K(V_K,\Theta_K)\cong \Cent_k(V,\Theta)\otimes_kK, $ where $ K $ is a field extension of $k.$ 

\medskip

For many $d$-linear spaces the center reduces just to the ground field $k.$ In fact, one can prove that if  $(V,\Theta)$ is 2-regular, then
$\Cent_k(\Theta) = k $, see \cite{kaw}.

\medskip

We mention another $d$-linear space that has center equal to $k$.
\begin{ex}\label{exscsatrace} Let ${\mathcal A}$ be a central simple algebra over $k,$ and let $\trace$ denote its reduced trace. 
Let $({\mathcal A}, T^d)$ be the $d$-linear space corresponding to the degree $d$ form  $\trace(x^d)$ on ${\mathcal A}.$ Explicitly, $T^d$ is given by 
$T^d(x_1,\ldots,x_d):=\frac{1}{d!}\trace(\sum_{\pi\in\mathfrak{S}_d}x_{\pi(1)}\cdots x_{\pi(d)})$ where $\mathfrak{S}_d$ denotes the symmetric group of degree $d.$ 
Then $\Cent_k({\mathcal A},T^d) \cong k,$ see  \cite[Theorem 2.2]{os1}.
\end{ex}

We now give another $d$-linear space that is constructed in quite a similar way to the space of the previous example and still has a very different center. Instead of 
considering a central simple algebra over $k$ we this time
consider $K,$ a finite separable field extension of $k.$ 
\begin{ex}\label{exseptrace}
Let $K$ be a separable finite extension of $k$ and let $b\in K, \, b\ne 0.$ Define the d-linear map $\Psi_b:K\times\cdots\times K\longrightarrow k$ by $\Psi_b(x_1,\ldots,x_d):= \text{tr}_{K/k}(bx_1\cdots x_d).$ The space $(K,\Psi_b)$ is regular, indecomposable and has center isomorphic to $K,$ see \cite{ha,hapa}. 
\end{ex}

Hence in general, the center of a regular $d$-linear space may be small as in example \ref{exscsatrace} or large as in example \ref{exseptrace}. 
Indeed, the center of the latter example is as large as possible, that is, it is a maximal commutative subalgebra of $\End_k(V), $ see 
Lemma {\ref{maximality}} below.
The forms that we introduce in the next section as our main object of study also have center with this maximality property, so let us  
formalize it:

\begin{definition}
Let  $(V,\Theta)$ be a regular  $d$-linear space over $k.$ We say that $(V,\Theta)$ has maximal center if $\Cent_k(V,\Theta)$ is a maximal commutative algebra in $\End_k(V).$ 
\end{definition}
\begin{lem}{\label{maximality}} 
1) The orthogonal sum 
of $d$-linear spaces with maximal center
has maximal center. Maximality is preserved under scalar extension.

2) Let $(K,\Psi_b)$ be as in Example \ref{exseptrace}. Then $(K,\Psi_b)$ has maximal center.

3) Let $K$ be finite algebraic field extension of $k$ and $(V,\Theta)$ a regular indecomposable d-linear space over $K$ with maximal center.
Let $s:K\longrightarrow k$ be a non-zero $k$-linear map and let be the $d$-linear space over $k$ with map $ (v_1, \ldots, v_d ) \mapsto 
s (\Theta(v_1, \ldots, v_d )) $. Then $(V,s\Theta)$ is regular, indecomposable and has maximal center. Moreover $\Cent_k(V,s\Theta) =\Cent_K(V,\Theta)$.

4) Let $({\mathcal A}, T^d)$ be as in Example  \ref{exscsatrace}, then its center is not maximal.
\end{lem}
\begin{proof}
1) A straightforward computation.

2) Identify $K$ with $\{m_a : a\in K\}\subset \End_k(K)$ where $m_a$ denotes multiplication by $a.$ Suppose $f\in \End_k(K)$ is such that  $f\circ m_a = m_a\circ f$ for all $a\in K.$ Then for each $a\in K$ we have  $\Psi_b(f(a),x_2\ldots,x_d)=\Psi_b(f(m_a(1)),x_2\ldots,x_d)=\Psi_b(m_a(f(1)),x_2\ldots,x_d) = \Psi_b(af(1),x_2\ldots,x_d).$ Since $\Psi_b$ is regular this implies that $f(a)=af(1)$ for each $a\in K,$ with $f(1)\in K$ and hence $f=m_{f(1)}\in K$ as needed.
%
%

3) The $d$-linear space $(V,s\Theta)$ over $k$ is regular and indecomposable by \cite[Proposition 3]{pu} and has $\Cent_k(V,s\Theta) =\Cent_K(V,\Theta)$ 
by \cite[Lemma 4.2 v)]{rup}. Since $K\subset \Cent_K(V,\Theta)$  one sees that any $k$-linear map on $V$ which commutes with every element in $\Cent_K(V,\Theta)$  is also 
$K$-linear and so $\Cent_k(V,s\Theta)$ is maximal.

4) The center is not maximal since scalar multiplication by 
an element of $k$ commutes with all $\End_k({\mathcal A}).$
\end{proof}

For more details on the forms $(V,s\Theta)$ given in 3) one may consult \cite[Lemma 2.7]{ha},  
\cite[Definition 2.3 iv)]{rup} 
and \cite[Proposition 3]{pu}.

\medskip


The {\em orthogonal} or {\em automorphism group} of a (regular) $d$-linear space $(V,\Theta)$ 
is defined as the set of $k$-linear bijections of $ V $ that leave $\Theta$ invariant: 
\begin{definition}
Let $(V,\Theta)$ be a (regular) d-linear space over $k$. The orthogonal group of $(V,\Theta)$ is 
\[\mathcal{O}(\Theta)=\{\sigma\in \GL_k(V)\,|\, \Theta(\sigma(v_1),\ldots,\sigma(v_d))=
\Theta(v_1,\ldots,v_d) \text{ for all } v_i\in V\}.\]
\end{definition}
\begin{ex}
The orthogonal group of the form given in example \ref{exscsatrace} is computed in \cite[Theorem 3.1]{os1} and is infinite if $k$ is. On the other 
hand the orthogonal group of example \ref{exseptrace} is finite, see \cite[Theorem 3.12]{hapa}.
\end{ex}

\medskip

Using the correspondence between $d$-linear spaces 
and homogeneous polynomials given in the beginning of this section, one sees 
that the orthogonal group $\mathcal{O}(\Theta)$ of $(V,\Theta)$ is a linear algebraic group. Hence there is a 
Lie algebra associated with $\mathcal{O}(\Theta)$ that we shall denote $\liet.$ By \cite[Proposition 4.2]{keet}, $\liet $ can 
be defined/described directly as follows: 
\begin{definition}\label{def-lie-algebra}
The Lie algebra  $\liet$ associated with $ (V,\Theta) $ is 
a subalgebra of $\mathfrak{gl}(V)$. It can be described as 
\[\liet=\{L\in \mathfrak{gl}(V)\,| \sum_{i=1}^d\Theta(v_1,\ldots,L(v_i),\ldots,v_d)= 0 \text{ for all } v_1,\ldots,v_d\in V\}\]
\end{definition}

\begin{ex}\label{lieexscsatrace}
The orthogonal groups of the spaces defined in Example \ref{exscsatrace} are computed in \cite{os1} from which one may derive their Lie algebras. Let us use 
the above description to show that they are at least nontrivial.
For $a\in {\mathcal A}$, let $L_a$ denote the left multiplication by $ a$ in $ {\mathcal A} $ and 
$R_a$ denote the right multiplication by $ a$ in $ {\mathcal A}. $ 
One checks that $\sum_{i=1}^dT^d(x_1,\ldots,L_a(x_i),\ldots,x_d) = \sum_{i=1}^dT^d(x_1,\ldots,R_a(x_i),\ldots,x_d),$ and hence by Definition 
\ref{def-lie-algebra} we have $\ad_a=R_a-L_a\in \lie_{T^d}.$ Then the 
Lie algebra has dimension at least $\dim_k{\mathcal A} -1.$
\end{ex}

On the other hand, by combining the next Proposition and example, one gets that the Lie algebras of the spaces in example 
\ref{exseptrace} are trivial.

\begin{prop}\label{dir-sum-ext}
Let $(V,\Theta)$ be a regular $d$-linear space (with $d>2$ as always) over $k$ and  $\liet$ be its Lie algebra.
\begin{enumerate}
\item Suppose that $(V,\Theta)\cong (V_1,\Theta_1)\perp (V_2,\Theta_2)$. Then $\liet\cong \lie_{\Theta_1}\times\lie_{\Theta_2}.$
\item For a field extension $ K/k$, one has $\lie_{\Theta_K}\cong \liet\otimes_kK.$
\end{enumerate}
\end{prop}
\begin{proof}
(1) See \cite[Proposition 4.3]{keet}.
(2) follows easily from Definition \ref{def-lie-algebra}.
\end{proof}  
\begin{ex}\label{lie-diagonal}
Suppose $\Theta$ is the d-linear form corresponding to the
homogeneous polynomial $F=a_1x_1^d+\cdots+a_nx_n^d$ of degree $d\geq 3 $ over $k$. Let $\{e_1,\ldots,e_n\}$ be the the canonical basis of $V=k^n.$ 
Then we have that $\Theta(e_{i_1},\ldots,e_{i_d})= a_{i_1}$ if $i_1=i_2=\cdots=i_d$ and $0$ otherwise.
Thus we have an  orthogonal decomposition $V=ke_1\perp ke_2\perp\ldots\perp ke_n.$ Proposition \ref{dir-sum-ext} now implies that $\liet\cong \lie_{\Theta_1}\times\lie_{\Theta_2}\times\cdots\times\lie_{\Theta_n},$ where $\Theta_i$
denotes the restriction of $\Theta$ to $ke_i.$ Since $ \Char k > d $ or $ \Char k = 0 $ we get $\lie_{\Theta_i} =0 $ and 
hence $ \liet =0 $. 
\end{ex}
A key property of the center of a $d$-linear space that we shall need in the following is its Lie module structure. 

\begin{prop}\label{action-lie-center}
Suppose $(V,\Theta)$ is a $d$-linear space over $k$ with center $\Cent_k(V,\Theta)$.
For $f\in \lie_{\Theta}$ and $\phi\in \Cent_k(\Theta)$ we have
$[f,\phi]\in  \Cent_k(V,\Theta)$ where the bracket is the commutator. This induces a 
module structure for $\lie_{\Theta}$ on $\Cent_k(\Theta).$ 
\end{prop}
\begin{proof}
The proof of $[f,\phi]\in  \Cent_k(V,\Theta)$ is a direct calculation involving 
the definitions/descriptions of $\Cent_k(V,\Theta)$ and $\lie_{\Theta}$. The module structure on $\Cent_k(V,\Theta)$
follows then from the basic fact that the commutator satisfies the Jacobi identity.
\end{proof}  

We now describe how the elements of the center can be used to twist the Lie algebra structure on $\lie_{\Theta}.$
\begin{prop}\label{gen-lie-center}
Let $ \psi \in \Cent_k(\Theta) $. Then there is a new Lie algebra structure on
$ \lie_{\Theta} $ given by $$ [f,g]_{\psi}  := f \psi g - g \psi f \,\,\,\,\,\,\,\,\,\,\,  \mbox{for} \,\,\,\, f,g \in \lie_{\Theta} $$
We denote $\lie_{\Theta}$ with this Lie algebra structure by $ \lie^{\psi} = \lie_{\Theta}^{\psi} $.
\end{prop}

\begin{proof}

It is clear that $ [f,g]_{\psi} $ is  bilinear and satisfies $ [f,f]_{\psi} = 0  $ whereas
the Jacobi identity is a simple calculation based on the definition.
Let us therefore verify that $ [f,g]_{\psi} \in  \lie_{\Theta} $ for
any $ f,g \in \lie_{\Theta} $. We must show that
\begin{equation}\label{longlie}
\Theta([f,g]_{\psi}v_1,v_2,\ldots, v_d) +
\Theta(v_1, [f,g]_{\psi}v_2,\ldots, v_d) + \ldots +
\Theta(v_1, \ldots , [f,g]_{\psi} v_d ) = 0
\end{equation}
For the first term of the above summation we have 
\begin{equation*}
\begin{aligned}
\Theta([f,g]_{\psi}v_1,v_2,\ldots, v_d) =
\Theta( (f \psi g - g \psi f ) v_1,v_2,\ldots, v_d) = \\
- \Theta( \psi g  v_1,f v_2,\ldots, v_d) - \ldots - \Theta( \psi g  v_1, v_2,\ldots, f v_d) +  \\
\Theta(  \psi f   v_1,g v_2,\ldots, v_d)  +  \ldots +  \Theta(  \psi f   v_1,v_2,\ldots,g v_d) = \\
\Theta(  f   v_1,g v_2,\ldots, \psi v_d)  +  \ldots +  \Theta(  f v_1, \psi v_2,\ldots,g v_d) \\
- \Theta( g v_1,f v_2,\ldots, \psi v_d) - \ldots - \Theta(  g  v_1, \psi v_2,\ldots, f v_d) \\
\end{aligned}
\end{equation*}
where in the last equality we used that $ \psi \in \Cent_k(\Theta) $ to move it to an `unoccupied' slot.
Expanding in the same way the other terms of (\ref{longlie})  we obtain similar expressions. 
Therefore, for each pair of distinct indices $ (i,j) $, there are exactly two 
terms in (\ref{longlie}) of the form $ \Theta(v_1, \ldots , f v_i, \ldots ,g v_j, \ldots , \psi v_l , \ldots,  v_d ) $ 
for $ l \not= i, j $, 
one term from the expansion of 
$ \Theta(v_1, \ldots , [f,g]_{\psi} v_i, \ldots  ,v_j, \ldots , v_d ) $ and another 
from the expansion of 
$ \Theta(v_1, \ldots ,  v_i, \ldots  ,[f,g]_{\psi} v_j, \ldots , v_d ).$ Since
the two terms have opposite signs, their sum is zero. 
\end{proof}

\section{Cyclic spaces}\label{cycsec}
We shall from now on focus on the ``cyclic spaces,'' 
which are the $d$-linear spaces whose center contains  a cyclic map.
In this section we recall the definition and basic facts of these spaces and then go on to calculate their Lie algebras, using 
the results from the previous section.
Somewhat surprisingly, we find that they are closely related to the Witt algebra.

\medskip
We assume from now that $ k $ is algebraically closed and $\Char k =0 $ even though
some of the results may hold in a greater generality. 

\medskip
\begin{definition}
Let $(V, \Theta) $ be a $d$-linear space. We say that $(V, \Theta) $ is cyclic if 
$ \Cent_k(\Theta) $ contains a cyclic element $ \psi $ for $ \EEnd(V) .$ This means that
there exists $v \in V$ such that $ V   = \spn \{ v, \psi(v), \psi^2(v), \ldots, \psi^N(v) \} $ 
for some $ N.$

\end{definition}

We need the following result from \cite{os}:
\begin{prop}\label{reich} Let $(V, \Theta) $ be regular and cyclic 
with cyclic element $\psi \in \Cent_k(\Theta).$ Then $\Cent_k(\Theta)= k[\psi]$ and is of  
dimension $n$ (where $ n = \dim V $).
\end{prop}


Regular indecomposable cyclic spaces exist in any dimension $n\geq 2$ and for any degree $d>2.$ 
Moreover for fixed $n,d$ they are unique up to multiplication by a scalar. 
There is a concrete construction of them, due to Reichstein \cite{rich}, which we shall 
recall now: 
\begin{definition}
\label{defrich}
Let $\{ v_1, \ldots , v_n \}$ be a basis of the vector space $V_n$ and $d>2$ an integer. 
Let $\Theta_d$ be the $d$-linear form defined by:
$$
\Theta_d(v_{i_1}, v_{i_2},\ldots, v_{i_d} ) = \left\{
\begin{array}{ll}
1 & \mbox{if} \,\,\, i_1 + i_2 + \ldots + i_d = (d-1)n +1 \\
0 & \mbox{otherwise}
\end{array}
\right.
$$
Then $ (V_n, \Theta_d ) $ is regular, indecomposable and cyclic. Indeed, 
let $\psi:V_n\longrightarrow V_n$ be the linear map defined by $ \psi(v_i) = v_{i-1} $ (where 
we define $ v_0 := 0 $).
Then  $\psi $ belongs to $ \Cent_k(V_n,\Theta_d)$ and is cyclic. 
Denote the center of $(V_n,\Theta_d)$ by $\Cent_k(n,d)$ 
and the Lie algebra by $\liend.$
\end{definition}


We know by Proposition \ref{reich} that $\Cent_k(n,d) = k[\psi]. $
The next series of results aim at computing $\liend.$ 

\begin{lem}\label{max-center-cyclic}
$ (V_n, \Theta_d ) $ has maximal center. 
\end{lem}
\begin{proof}
Let $ \psi \in \Cent_k(n,d) $ be as above and 
consider $ V_n $ as $ k[X] $-module through $ X \mapsto \psi $.
Since $ \psi $ is cyclic, this gives an isomorphism of $ k[X] $-modules  
$ V_n \cong k[X]/(X^n) \cong  \spn_k \{1,X, \ldots , X^{n-1} \} $. 
Hence if $ f \in \End_k (V_n) $ commutes with $ \psi $, it may be viewed as a linear map on 
$ \spn_k \{1,X, \ldots , X^{n-1} \} $ commuting with $ X$. We then have $ f(1) = P(X) $ for some $ P(X) \in k[X] $
and $ f(X^i)  = X^i f(1) = P(X) X^i   \mod X^n .$ In other words, $ f $ acts on $ V_n $ as 
multiplication by $ P(\psi) $, that is $ f \in \Cent_k(n,d) $ as claimed. 
\end{proof}
 
\begin{prop}\label{solvab-lie-cyclic}
Assume that $(V_n,\Theta_d)$ is as above. Then $\liend $ 
is a subalgebra of the upper triangular matrices of $\mathfrak{gl}_n(k)$ and hence solvable.
\end{prop}
\begin{proof}
Let $\psi$ and $\{v_1,\ldots,v_n\}$ be as in the above definition and take $f\in \liend.$ By Proposition \ref{action-lie-center} we have that $ [f, \psi]=  f\psi-\psi f\in \Cent_k(n,d).$  
But $[f, \psi ]$ 
is traceless and so it belongs to $ \psi \Cent_k(n,d) = \psi k[ \psi ] $.
Since $ [f , \psi^i ] =  \psi^{i-1} [f , \psi ] + [f , \psi^{i-1} ]  \psi $ for $ i \geq 1 $
we get now by induction that $ [f , \psi^i ] $ belongs to $ \psi^{i} k[ \psi ] $.
Set $ v:= v_n $ so that $ v_i = \psi^{n-i} v $. We then have 
$$ f v_i = f \psi^{n-i} v = \psi^{n-i} f v + [f, \psi^{n-i}] v \in \psi^{n-i} k[ \psi ] v 
= \spn \{ v_i, v_{i-1}, \ldots, v_1 \} $$
and the Proposition is proved.
\end{proof}

\begin{lem}\label{dim-lie-cyclic}
For any $d>2$ and $n\geq2$ we have 
\[\dim_k \liend< \dim_k \Cent_k(n,d) =n.\]
\end{lem}
\begin{proof}
We need only check the first inequality. 
Proposition 
\ref{action-lie-center} gives rise to a linear map  
$ i: \liend \rightarrow \Cent_k(n,d), \, \, 
f \mapsto  [f,\psi] $, where $ \psi $ is as above.
We show that $ i $ is injective. 

\medskip
Suppose that $f \in \liend \cap \ker i $. Then  
$f \in \Cent_k(n,d)$, since $ (V_n, \Theta_d ) $ has maximal center, 
that is $\Theta_d(f(u_1),u_2,\ldots,u_d)= \Theta_d(u_1,f(u_2),\ldots,u_d) = \ldots = 
\Theta_d(u_1,u_2,\ldots,f (u_d)) $ for all $u_1,\ldots,u_d\in V_n . $
But we also have 
$$\Theta_d(f(u_1),u_2,\ldots,u_d)+\Theta_d(u_1,f(u_2),\ldots,u_d)+\cdots
+\Theta_d(u_1,u_2,\ldots,f(u_d))=0$$
since $ f \in \liend $. 
Combining, we get $\Theta_d(f(u_1),u_2,\ldots,u_d)=0$ for all
$u_1,\ldots,u_d\in V_n$ (recall $ \Char k = 0 $). This implies that $f =  0 $ because $ (V_n, \Theta_d) $ is regular, and so indeed $ i $ is injective. Moreover, since $ \trace [f, \psi] = 0 $ we have $ 1 \notin \im i $ and the Lemma is proved. 
\end{proof}
For small $n$ and $d$ we can compute the Lie algebra 
$ \liend $ explicitly. 
\begin{ex}
We have
$\mathcal{L}(3,3)=\left \{
 \left[
\begin{smallmatrix}
 - 2a & b & 0 \\
0 & a & -2b \\
0 &  0 & 4a
\end{smallmatrix}
 \right]:a,b\in k\right \}$ and  
$\mathcal{L}(4,3)=\left \{
 \left[
\begin{smallmatrix}
a & -2b & c & 0 \\
0 & 0 & b & -2c \\
0 & 0 &  - a & 4b \\
0 & 0 & 0 &  - 2a
\end{smallmatrix}
 \right] : a,b,c\in k\right \}$
\end{ex}
These examples suggest that there might be a sort of 
inclusion of the $n$ dimensional case into the $n+1$ dimensional case. 
We shall show that this is indeed the case. 

\medskip
Let $\psi \in \Cent_k(n,d) $ be as above. Define the linear maps 
$$ \iota: V_n \rightarrow V_{n+1}, \,\,\,\, \iota(v_i) = v_i,  \,\,\,\,\,\,\,\,\,
\pi: V_{n+1} \rightarrow V_{n},  \,\,\,\,\,\, \pi(v_i) = v_{i-1}. $$ 
Then we have 
\begin{prop}\label{embedding}
The map
$ \rho:  \liend^{\psi} \rightarrow \lie(n+1,d) $
defined by $ f \mapsto \iota \circ f \circ \pi $ is a Lie algebra embedding, where 
$ \liend^{\psi} $ is the twisted Lie algebra structure on $ \liend $ introduced in 
Proposition \ref{gen-lie-center}.
\end{prop}
\begin{proof}
Certainly $ \rho $ is an injection of vector spaces. Since $  \pi \circ \iota = \psi $
we have that
$$
 [\rho(f),\rho(g)] = \iota f \pi \iota g \pi - \iota g \pi \iota f \pi =
\iota f \psi g \pi - \iota g \psi f  \pi  = \rho( [f,g]_{\psi})
$$
and hence it is a Lie algebra homomorphism.
It only remains to be proved that $ \rho(f) \in \lie(n+1,d) $.
But this follows from the following formula and its permutations: 
$$ \Theta_d( \iota v_{i_1}, v_{i_2}, \ldots , v_{i_d} ) =
\Theta_d( v_{i_1}, \pi v_{i_2},   \ldots ,\pi  v_{i_d} ) $$ 
for $ v_{i_1} \in
V_n$, $  v_{i_2}, \ldots , v_{i_d} \in V_{n+1} $. It can be read off from
the standard forms given in Definition \ref{defrich}. \end{proof}

\begin{lem}\label{diagonal-element}
Define $ D_n \in \EEnd_k(V_n)$ by $ D_n: v_{i} \mapsto  ( n - 1 - d(n-i) )
v_{i}, $ for all $i.$ Then $ D_n $ is a semisimple element of $ \liend.$
\end{lem}
\begin{proof}
By definition $ \Theta_d( v_{i_1}, v_{i_2},\ldots , v_{i_d}
) \not=  0 \Rightarrow \sum_j i_j  = (d-1)n + 1.$ But
\begin{eqnarray*}
\Theta_d(D_n v_{i_1}, v_{i_2},\ldots , v_{i_d} ) + 
\cdots +
\Theta_d(  v_{i_1}, v_{i_2},\ldots , D_n v_{i_d} ) = \\
d ( (n-1) -   \sum  (n- i_j ))  \, \Theta_d( v_{i_1}, v_{i_2},\ldots , v_{i_d} )  = \\
d (\, (n-1) -   d n +  \sum i_j \,  )  \, \Theta_d( v_{i_1}, v_{i_2},\ldots , v_{i_d} ) = 0,  
\end{eqnarray*}
and hence $ D_n \in
\liend $.
Moreover, since 
$ \liend $ is a linear Lie algebra and 
$ D_n $ acts diagonally in $ V_n $, it follows that it is a semisimple element
in the sense of Lie algebras. The Lemma is proved.
\end{proof}
\begin{lem}\label{diagonal-element-comm}
The following commutation rule holds in $ \liend $. $$ [D_n,\psi] = d \psi $$
\end{lem}
\begin{proof}
Apply $ D_n$ to the equation $ \psi v_i = v_{i-1} $.
\end{proof}
\begin{thm}\label{theorem-conjecture}
$ \lie(n+1,d) = \rho(\liend)\, \bigoplus \,kD_{n+1} $ and $ \dim_k \lie(n,d) = n-1.$ 
\end{thm}
\begin{proof}
We proceed by induction. Take 
$ f \in \lie(2,d) $ and write $f=\left [ \begin{smallmatrix} a_1 & a_2 \\ a_3 & a_4 \end{smallmatrix}\right ]$ 
%
with respect to the basis $\{ v_1, v_2 \}$ of $V_2.$ The condition 
$$ 
\begin{array}{rl}
\Theta_d(f(v_2),v_2,\ldots,v_2, v_1, v_1 ) +
\Theta_d(v_2,f(v_2),\ldots,v_2, v_1, v_1 ) &+ \\
\Theta_d(v_2,v_2,\ldots,f(v_1), v_1 )+  \Theta_d(v_2,v_2,\ldots,v_1, f(v_1) )&=0 
\end{array}
$$
implies that $ a_3= 0 $. Similarly one gets that $a_2=0 $ by acting on $ (v_2, v_2 \ldots , v_2 ) $ and finally 
$ a_1 + (d-1)a_4 = 0 $ by acting on $ (v_2, v_2 \ldots , v_2, v_1  ) $. We conclude that 
$  \lie(2,d) $ has dimension 1.

Assume now inductively that the Theorem is true for $ \lie(r,d)$ for $ r
\leq n $.  We then get by Proposition ~\ref{embedding} that $
\rho(\liend)$ consists of nilpotent matrices. Thus the vector
space sum of the Theorem is direct and so $ \dim
\lie(n+1,d) \geq \dim  \rho(\liend) + 1 $. The assertion
on the dimension now follows using Proposition
~\ref{dim-lie-cyclic} and from this we get that $ \rho(\liend) $
and $ D_{n+1}$ span $ \lie(n+1,d) $.
\end{proof}
\begin{rem}
Recall that the Witt algebra $ W$ is the Lie-algebra 
on generators $ \{ L_n \, | \, n \in \mathbb Z \} $ 
and relations $[ L_n, L_m ] = (m-n) L_{n+m}$.
\end{rem}

The following Theorem states that $ \liend $ is a quotient of a certain subalgebra of $ W$, namely the one
given by $ \{ L_n \, |\, n = 0, 1 , 2, \ldots \,   \} $. It is the main Theorem of this section.


\begin{thm}\label{Witt-algebra} 
Let $ \Dn := \frac{1}{d} D_n $ and $ X_r = (\Dn + r \frac{d-1}{d} I) \psi^r.$ Then we have

i) $ \{ X_r \, | \, r= 0, \ldots , n-2 \} $ is a basis of $ \liend $   \newline
ii) $ [X_r, X_s] = \left\{ \begin{array}{ll} (s-r)X_{s+r} & \mbox{if    }\,\,\,\, s+r < n-1 \\ 0 & \mbox{otherwise.} \end{array} \right.$
\end{thm}

\begin{proof} i). We first claim that $ X_r \not= 0 $ for $ r= 0, \ldots, n-2 $. But this is a consequence of the formula
$ X_r v_n = \frac{n-1-r}{d} v_{n-r} $
that follows easily from the definitions. Since $ \{v_i\, |\, i = 1,\ldots , n \}$ is a basis 
of $ V_n $ we even see from this 
that $\{ X_r \, |\,  r = 0,\ldots , n-2 \}$ is a linearly 
independent subset of $ \End_k (V_n) $.
So we just have to show 
$ X_r \in \lie(n,d) $. For this it is enough by linearity to check that
$$ \Theta_d(X_r v_{i_1}, v_{i_2},\ldots , v_{i_d} ) + \Theta_d(v_{i_1}, X_r
v_{i_2},\ldots , v_{i_d} )+  \cdots + \Theta_d( v_{i_1},\ldots ,X_r v_{i_d} ) = 0.  $$ 
By definition of $ \Theta_d $ we only need check the case 
$ i_1 + i_2 + \ldots + i_d = (d-1)n +1 +r.$

Now we have for all $ j = 1, \ldots , d \, $:
$$ \Theta_d( v_{i_1}, \ldots,  X_r v_{i_j},\ldots , v_{i_d} ) =  \frac{n - 1 }{d}   - (n-i_j+r) + r \frac{d-1}{d}.   $$
Summing up we find
$$ 
\begin{array}{r}
\Theta_d(X_r v_{i_1}, v_{i_2},\ldots , v_{i_d} ) +  \ldots + \Theta_d( v_{i_1}, v_{i_2},\ldots ,X_r v_{i_d} ) = \\
n-1 -nd -rd +r(d-1) + (d-1)n +1+r   =0 
\end{array}
$$
as claimed. 
\newline
ii). 
One proves first that $[ \Dn {\psi},  {\psi}^s ] = s  {\psi}^{s} $ by induction on $ s $ using 
the formula $[ \Dn  ,  {\psi} ] = {\psi} $ from  Lemma \ref{diagonal-element-comm} together with 
$[ \Dn  ,  {\psi}^s ] = \psi [ \Dn  ,  {\psi}^{s-1} ] + [ \Dn  ,  {\psi} ]{\psi}^{s-1} .$ 
From this one gets 
$[ \Dn {\psi}^r ,  {\psi}^s ] = s  {\psi}^{r+s} $ by induction on $ r$ and then 
$[ \Dn {\psi}^r , \Dn {\psi}^s ] = (s-r) \Dn {\psi}^{r+s}. $
This formula is valid 
for all integers $r, s \geq 0 .$ But now 
$$ 
\begin{array}{r}
[ X_r, X_s ] = [ (\Dn + r \frac{d-1}{d} I) \psi^r, (\Dn + s \frac{d-1}{d} I) \psi^s ] = \\

[\Dn {\psi}^r, \Dn {\psi}^s] +  s \frac{d-1}{d} [\Dn {\psi}^r ,  {\psi}^s ] - r \frac{d-1}{d} [\Dn {\psi}^s ,  {\psi}^r ]=  \\

(s-r) \Dn {\psi}^{r+s}+  (s^2 -r^2) \frac{d-1}{d} \psi^{s+r} = (s-r) X_{s+r} 
\end{array}
$$ 
as claimed.
Hence to get the expression of ii), we must check that $ X_r = 0 $ for $ r \geq n-1 $. 
For $ r \geq n $ this is clear since $ \psi^n =0 $. And for $ r = n-1 $ it follows from $ v_l =0 $ for $ l \leq 0 $ and from the formula 
$ X_r v_n = \frac{n-1-r}{d} v_{n-r} $ that we used in the proof of i). 
\end{proof}

\begin{rem}
Note that neither $\lie(n,d)$
nor $\Cent_k(n,d) $ depends on $ d.$
\end{rem}

\begin{rem}
One may ask if it is possible to realize the Witt algebra itself as the Lie algebra of a
higher degree form. Allowing infinite dimensional vector spaces, one possible way of doing so is
to let $ V $ 
be the vector space with 
basis $ \{ \, v_i \, |\,  i \in {\mathbb Z}\,  \} $ and define $ \Theta_d $ by 
$$
\Theta_d(v_{i_1}, v_{i_2},\ldots, v_{i_d} ) = \left\{
\begin{array}{ll}
1 & \mbox{if} \,\,\, i_1 + i_2 + \ldots + i_d = 0 \\
0 & \mbox{otherwise}
\end{array}
\right.
$$
We leave out the details.
\end{rem}

\begin{qu}
Is is possible to realize the Virasoro algebra, the central extension of the Witt algebra, 
as the Lie algebra of a higher degree form?
\end{qu}

\section{The Orthogonal Group}\label{autsec} 
In this section we use the results from the previous section to obtain a precise description 
of the orthogonal group in the cyclic case. 

\medskip
In general, in order to determine the orthogonal group of a regular space $ (V, \Theta) $ with $d>2$ one may assume 
that it is indecomposable. Indeed, 
under the action of an automorphism of $ (V, \Theta) $, 
the indecomposable components are just being permuted, 
as one sees by the uniqueness of the components.

\medskip
The following Lemma is the key point for the results of this section.
It is the group theoretical version of 
Proposition \ref{action-lie-center}.

\begin{lem}
Let $(V,\Theta)$ be a regular  $d$-linear space over $k.$ Then $\mathcal{O}(\Theta)$ acts on $\Cent_k(\Theta)$ by conjugation, that is for each $\sigma \in  \mathcal{O}(\Theta)$ and $f\in  \Cent_k(\Theta),$ we have $\, \sigma f\sigma^{-1}\in \Cent_k(\Theta)$. 
\end{lem}
\begin{proof}
Follows directly from the definitions.
\end{proof}

Let $\mu_d :=\{\zeta\in k \,|\, \zeta^d=1\},$ which by the assumptions on $ k $ 
has order $ d$. Note that $ \mu_d $ may 
be identified with a subgroup of $\mathcal{O}(\Theta).$  Let $\mathbf{G}$ be the group of automorphisms of the $k$-algebra $\Cent_k(\Theta).$ It is an algebraic group and 
the action of $\mathcal{O}(\Theta) $ on $ \Cent_k(\Theta) $ induces 
a homomorphism of algebraic groups $ \chi: \mathcal{O}(\Theta) \rightarrow \mathbf{G}. $ We can now formulate our next Theorem:



\begin{thm}
Let $(V,\Theta)$ be a regular indecomposable $d$-linear space over $k$ with maximal center.  Then 
$ \chi $ induces the following exact sequence of groups 
\[1\longrightarrow\mu_d\longrightarrow \mathcal{O}(\Theta) \stackrel{\chi}{\longrightarrow}\mathbf{G}\]
\end{thm}

\begin{proof}
Suppose that $\sigma\in \mathcal{O}(\Theta)$ satisfies $\chi (\sigma)=1,$ that is  
$\sigma\rho\sigma^{-1}=\rho$ for all $\rho \in  \Cent_k(\Theta).$ Then $\sigma$ 
is in the centralizer of $\Cent_k(\Theta),$ which by maximality implies that 
$\sigma \in \Cent_k(\Theta).$
Hence for all $v_1,\ldots,v_d\in V$, we have
$$\Theta(v_1,\ldots,v_d)=\Theta(\sigma(v_1),\ldots,\sigma(v_d))=\Theta(v_1,\ldots,\sigma^d(v_d)). $$
By regularity of the space we get from this  $\sigma^d=id_V.$ Therefore the minimum polynomial of $\sigma$ divides $X^d-1\ $ in $k[X]$ 
and so the eigenvalues of $\sigma$ are multiplicity free and $ \sigma $ is diagonalizable.
If $\sigma$ had more than one eigenvalue,  then $ (V, \Theta) $ 
would be decomposable by 
\cite[Lemma 2.6]{os}.
Hence $\sigma$ has exactly one eigenvalue and we conclude that 
$\sigma=\zeta id_V$ with $\zeta^d=1$ as needed.
\end{proof}

We now return to the cyclic spaces $(V_n,\Theta_d)$ of the previous section.
The following is our main Theorem. 

\begin{thm}\label{orthcyc}
Let $(V_n,\Theta_d)$ be the regular, indecomposable cyclic $d$-linear space over $k$ 
with cyclic element $ \psi \in \Cent_k(n,d)$. Then the homomorphism $ \chi $ induces the following short exact sequence of groups
\[1\longrightarrow\mu_d\longrightarrow \mathcal{O}(\Theta_d) \stackrel{\chi}\longrightarrow \mathbf{G} {\longrightarrow} 1\]
\end{thm}

Before we can give the proof of the Theorem we need a couple of preparatory lemmas. 
Let us first analyse the group $\mathbf{G}$ in more detail.
Consider the set $ \mathbf{G^{\prime}} := k^\times\times k^{n-2}$ and 
define $ \rho_1 : \mathbf{G^{\prime}} \rightarrow \Cent_k(n,d) $
by 
$  a = (a_1 , a_2 \ldots , a_{n-1})    \mapsto a_1 \psi + a_2 \psi ^2  + \ldots a_{n-1} \psi ^{n-1}. $
Define moreover $ \rho_2 : \mathbf{G^{\prime}} \rightarrow \Matnk $ by 
$$ \rho_2 (a) := \left[\begin{smallmatrix}
1                       & 0                     & 0                     & \cdots                & 0 \\
0                       & a_1           & 0             & \cdots                & 0     \\
0                       & a_2           & a_1^2 &                               & 0     \\\
\vdots  & \vdots        &       & \ddots    & \vdots\\    
0                       & a_n           & \star &\star                  & a_1^{n-1}     
\end{smallmatrix}\right] $$
where the $i$'th column consists of the  entries of $\rho_1 (a)^i  $. This induces an operation $ m: \mathbf{G^{\prime}} \times \mathbf{G^{\prime}} \rightarrow 
\mathbf{G^{\prime}} $
where $ m(a, b) $ is the second column, with the first entry deleted, 
of the matrix product $   \rho_2(a)  \rho_2(b) $ 
\begin{lem}
Let  $(V_n,\Theta_d)$ be as above.
Then $\mathbf{G^{\prime}} $ is a group with multiplication given by $ m $ and neutral element $ e := (1,0,0, \ldots, 0 )$. 
The groups $\mathbf{G^{\prime}} $ and $\mathbf{G} $ are isomorphic.
\end{lem}  
\begin{proof}
Notice first that $ a \in \mathbf{G} $ is completely determined by $a(\psi)$ since
$ a$ is a $k$-automorphism of $ \Cent_k(\Theta) .$
Moreover, 
$a(\psi)= a_1\psi+\cdots+a_{n-1}\psi^{n-1}$ for $a_i\in k \text{ and } a_1\neq 0$ 
because $ \psi^{n} = 0 $ and $ \psi^{j} \not= 0 $ for $ j \leq n-1 $.
Hence $a $ determines a vector $f(a) = (a_1,\ldots,a_{n-1}) \in
k^\times \times k^{n-2}$, and so we have a map $ f : \mathbf{G} \rightarrow  \mathbf{G^{\prime}} .  $
Since  
$ \Cent_k(\Theta) $ is generated by $ \psi$, we deduce that $ f $ is bijective. 
The multiplication and 
neutral element of $ \mathbf{G} $ 
give by transport of structure via $f$ exactly 
the multiplication $ m $ and neutral element $e $ on $ \mathbf{G^{\prime}} $ and the Lemma is
proved.
\end{proof}

\begin{lem}{\label{connected}}
$\mathbf{G}$ is a connected algebraic group.
\end{lem}
\begin{proof}


From the previous Lemma we know that $ \rho_2:  k^{\times} \times k^{n-2 } \rightarrow Mat_n(k) $ 
makes $ k^{\times} \times k^{n-2 }$ in bijection 
with its image $\mathbf{G} $. The inverse map is given by projection on the second column.
Hence, as a variety $\mathbf{G} $ is isomorphic to $ k^{\times} \times k^{n-2 }$ and
so $\mathbf{G} $ is connected.
\end{proof}

\begin{lem}
The dimension of $\mathbf{G}$ is $n-1.$ 
\end{lem}
\begin{proof}
By \cite[section 13.2]{hump} the Lie algebra $\mathfrak g$ of $\mathbf{G}$ is the Lie algebra of 
derivations of $k[\psi].$ 
But any derivation  $d:k[\psi]\longrightarrow k[\psi]$ is uniquely determined by $d(\psi)$ and so
$\mathfrak g$ has dimension $ n-1 $. Hence also $ \mathbf{G} $ has dimension $ n-1 $.
(Alternatively, one can argue directly using $ \mathbf{G} \cong k^{\times } \times k^{n-2} \, $).
\end{proof}
We are now able to prove Theorem \ref{orthcyc}.
\begin{proof}
We need only show $\chi:\mathcal{O}(\Theta_d)\longrightarrow \mathbf{G}$ is surjective. Set $\mathbf{H}:= \chi(\mathcal{O}(\Theta_d)).$ It is a closed subgroup 
of $\mathbf{G}$, \cite[Section 7.4, Proposition B]{hump}.
By Theorem \ref{theorem-conjecture}, the Lie algebra of $ (V_n, \Theta_d )$ has dimension $n-1$ 
and hence also 
$\mathcal{O}(\Theta_d)$ has dimension $ n-1$.  Then $\mathbf{G}/\mathbf{H}$ is a variety of 
dimension $0$ and so $\mathbf{H}$ contains the identity component $\mathbf{G}^0$, 
\cite[Proposition 7.3 (b)]{hump}.
But by  
Lemma \ref{connected}, $\mathbf{G}$ is connected and therefore $\mathbf{G}=\mathbf{H}.$ 
\end{proof}

By Theorem \ref{orthcyc}, $ \chi $ induces an isomorphism $\mathcal{O}(\Theta_d)/ \mu_d \cong \mathbf{G}  $. 
In general, however, we don't know how to describe the 
inverse map of $ \chi $ directly without passing through the Lie algebras. 

\medskip
Of course for small values of $n$ and $ d$, one can use an ad-hoc approach. 
Let's give an explicit example.  
Assume that $n=d=3.$ From the above description, $\mathbf{G}$ is generated by the elements 
$\rho = \left[\begin{smallmatrix}
                                        1 & 0 & 0\\
                                        0 & \alpha & 0\\
                                        0 & 0 & \alpha^2
\end{smallmatrix}\right] \text{ and } \epsilon 
 = \left[\begin{smallmatrix}
                                        1 & 0 & 0\\
                                        0 & 1 & 0\\
                                        0 & \beta & 1
\end{smallmatrix}\right].$ 
We look for 
$\sigma \in \mathcal{O}(\Theta_3)$ such that $\chi(\sigma)=\rho$ 
and 
$\tau \in \mathcal{O}(\Theta)$ such that $\chi(\sigma)=\epsilon.$
Let's start with $ \sigma$.
\medskip

Such $ \sigma $ must satisfy $\chi(\sigma)(\psi)=\sigma\psi\sigma^{-1}=\alpha\psi.$ Combining with $\psi(v_i)=v_{i-1}$ and writing $\sigma(v_3)=av_1+bv_2+cv_3$ an easy computation shows that 
$\sigma = \left[\begin{smallmatrix}
                                        \alpha^2c & \alpha b  & a\\
                                        0 & \alpha c& b\\
                                        0 & 0 & c
\end{smallmatrix}\right].$ 
Then, using that $\sigma$ is an automorphism of $\Theta_3$ we obtain the following system of equations: $\alpha^2 c^3=1,\, \alpha bc^2=0,\, b^2c+ac^2= 0.$ Solving this system and denoting  by $\sqrt[3]{\alpha}$ a fixed cubic root of $\alpha$ we see that the matrix of $\sigma$ is  
\[\sigma=\left[\begin{smallmatrix}
                                        \alpha\sqrt[3]{\alpha} & 0 & 0\\
                                        0 & \sqrt[3]{\alpha}& 0\\
                                        0 & 0 & \sqrt[3]{\alpha}/\alpha
\end{smallmatrix}\right]\]
To obtain $\tau \in \mathcal{O}(\Theta)$ we proceed 
similarly and find \[\tau=\left[\begin{smallmatrix}
                                        \zeta & 2\beta\zeta/3 & -\beta^2\zeta^2/9\\
                                        0 & \zeta &  -\beta\zeta/3\\
                                        0 & 0 & \zeta
\end{smallmatrix}\right],\]
where $\zeta$ a cubic root of 1.
%
%
%


\begin{thebibliography}{WWW}
\bibitem[C]{cl} A. Chlebowicz: Certain finite groups as automorphism groups of forms of higher degree, \textit{Linear Algebra Appl.} {\bf 419} (2006), 326-330.
\bibitem[CSW]{csw} A. Chlebowicz, A. S\l adek, A. Weso\l owski: Automorphisms of certain forms of higher degree over ordered fields, \textit{Linear Algebra Appl.} {\bf 331} (2001), 145-153.
\bibitem[H]{ha} D. Harrison: A Grothendieck ring of higher degree
forms, \textit{J. Algebra} {\bf 35} (1975), 123--138.
\bibitem[HP]{hapa} D. Harrison and B. Pareigis: Witt rings of higher degree
forms, \textit{Comm. Algebra} {\bf 16}(6) (1988), 1275-1313.
\bibitem[Hu]{hump} J. E. Humphreys: Linear Algebraic Groups, Springer-Verlag, GTM {\bf 21}, New York, 1998.
\bibitem[K]{keet} A. Keet: The Lie Algebra of a  Higher Degree Form and a Schur Functor, \textit{Comm. Algebra} {\bf 22}(5) (1994), 1577-1601.
\bibitem[KW]{kaw}T. Kanzaki and Y. Watanabe: Determinants of $r$-fold symmetric multilinear forms, \textit{J. Algebra} {\bf124} (1989),
219-229.
\bibitem[MM]{matsu} H. Matsumura, P. Monsky: On the automorphisms of hypersurfaces, \textit{J. Math. Kyoto Univ} {\bf 3} (1964), 347-361.
\bibitem[OS]{os} M. O'Ryan, D. B. Shapiro: Centers of higher degree forms, \textit{Linear Alg.  Appl.} {\bf 371} (2003), 301-314.
\bibitem [OS1]{os1}M. O'Ryan and D. Shapiro: On trace forms of Higher
Degree, \textit{Linear Alg. Appl.} {\bf 246} (1996), 313-333.
\bibitem [Pr]{pros}A. Pr\'{o}szy\'{n}ski: On orthogonal decomposition of homogeneous polynomials, \textit{Fund. Math.} {\bf 98} (1978), 201-217.
\bibitem [P]{pu} S. Pumpl\"{u}n: Indecomposable forms of higher degree,  \textit{Math. Zeit.} {\bf 253} N.2  (2006), 347-360.
\bibitem[R]{rich} B. Reichstein: On Waring's problem for cubic
forms, \textit{Linear Alg. Appl.} {\bf 160} (1992), 1-61.
\bibitem[Ru]{rup} C. Rupprecht: Cohomological invariants for higher degree forms, PhD Thesis, Universit\"{a}t Regensburg, 2003.

\bibitem[SC]{schnei} J. E. Schneider: Orthogonal groups  of nonsingular forms of higher degree, \textit{J. Algebra} {\bf 27} (1973), 112-116.
\bibitem[SW]{sla} A. S\l adek, A. Weso\l owski: Clifford-Littlewood-Eckmann groups as orthogonal groups of forms of higher degree, \textit{Annales mathematicae Silesianae} {\bf 12} (1998), 93-103.
\bibitem[S]{spr} T. A. Springer: Linear algebraic groups, Second Edition, Birkh\"auser, PM {\bf 9}, 1998
\bibitem[Sum]{summ} L. Summerer: Decomposable forms and automorphisms, \textit{J. Number Theory} {\bf 99} (2003), 232-254.
\bibitem[Suz]{suzu} H. Suzuki: Automorphism groups of multilinear maps, \textit{Osaka J. Math} {\bf 20} (1983), 659-673. 
\end{thebibliography}
\end{document}